\documentclass[runningheads,numbook]{svjour2}
\smartqed  % flush right qed marks, e.g. at end of proof

\usepackage{graphicx,epsfig}
\usepackage{amsmath,amscd,amsfonts}

\newtheorem{thm}{Theorem}[section]
\newtheorem{cor}[thm]{Corollary}
\newtheorem{lem}[thm]{Lemma}
\newtheorem{prop}[thm]{Proposition}

\numberwithin{equation}{section}
\newcommand{\CalT}{\mathcal {T}}
\newcommand{\CalS}{\mathcal {S}}
\newcommand{\CalB}{\mathcal {B}}
\newcommand{\CalM}{\mathcal {M}}
\newcommand{\lra}{\longrightarrow}
\journalname{Mathematische Annalen}
\begin{document}

\title{Braids, mapping class groups, and categorical delooping
\thanks{The first author was supported by Inha University research grant.}
}


\author{Yongjin Song  \and  Ulrike Tillmann}


\institute{Yongjin Song \at Department of Mathematics, Inha
University, 253 Yonghyun-dong, Nam-gu, Incheon 402-751, Korea \\
\email{yjsong@inha.ac.kr} \and Ulrike Tillmann \at Mathematical
Institute, 24-29 St Giles Street, Oxford OX1 3LB, UK \\
\email{tillmann@maths.ox.ac.uk}}

\date{Received: date / Revised: date}

\maketitle

\begin{abstract}
Dehn twists around simple closed curves in  oriented surfaces
satisfy the braid relations. This gives rise to a group theoretic
map $\phi: \beta_{2g} \to \Gamma _{g,1}$ from the braid group to
the mapping class group. We prove here that this map is trivial in
homology with any trivial coefficients in degrees less than $g/2$.
In particular this proves an old conjecture of J. Harer. The main
tool is categorical delooping in the spirit of \cite{T1}. By
extending the homomorphism to a functor of monoidal 2-categories,
$\phi$ is seen to induce a map of double loop spaces on the plus
construction of the classifying spaces. Any such map is
null-homotopic. In an appendix we show that geometrically defined
homomorphisms from the braid group to the mapping class group
behave similarly in stable homology. \subclass{55P48 \and 57M50
\and 55R37}
\end{abstract}

\section{Introduction}

Let  $\beta _k$ be Artin's braid group on $k$
strings \cite{A}. In its standard presentation $\beta _k$ has
generators $\sigma _1, \dots, \sigma _{k-1}$ subject to  the
following relations:
$$
\begin{aligned}
\sigma _i \sigma _j = \sigma _j \sigma _i \quad & \text {for } |i-j| \geq 2; \\
\sigma _i \sigma _{i+1} \sigma _i = \sigma _{i+1} \sigma _i \sigma
_{i+1} \quad & \text{for }  i= 1, \dots, k-2.
\end{aligned}
$$
Let $\Gamma _{g,1}$ be the mapping class group of an oriented
surface with one boundary component. Wajnryb \cite{W}, \cite{BW}
gave a presentation for $\Gamma _{g,1}$ with generators the  Dehn
twists $\alpha _i, i= 0, \dots ,2g-1,$ and $\delta$ around the
simple closed curves $a_i, i= 1, \dots , 2g-1$, and $d$ as
depicted in Figure 1. Among others, $\alpha _1 , \dots , \alpha
_{2g-1}$ satisfy the same relations as the braid group generators
above.
%Indeed, whenever two simple closed curves $a$ and $d$
%do not intersect, the associated
%Dehn twists $\alpha$
%and $\gamma$ commute, and if they intersect once, the braid
%relation $\alpha \gamma \alpha = \gamma \alpha \gamma$ is satisfied.
\begin{center}
\epsfig{figure=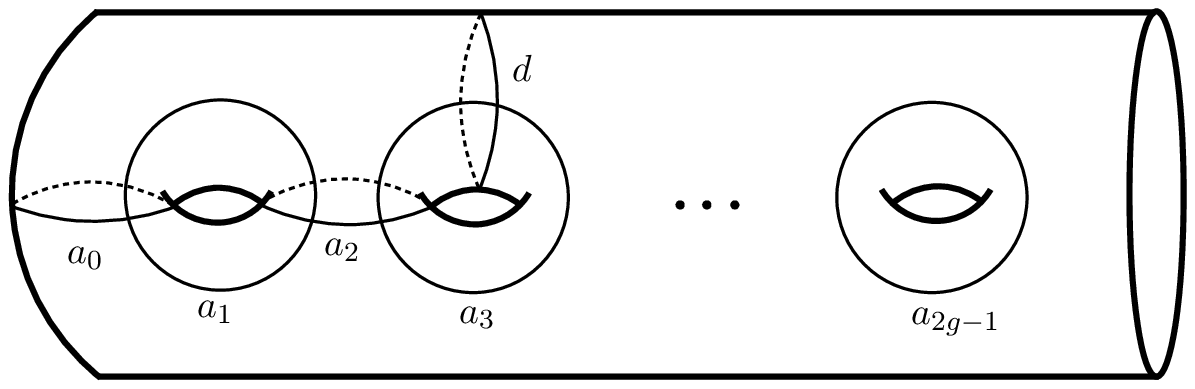, width=0.7\linewidth}

Figure 1.
\end{center}

\vskip .1in Hence, the assignment $\sigma _i \mapsto \alpha _i$
defines a purely group-theoretic homomorphism
$$
\phi: \beta _{2g} \longrightarrow  \Gamma _{g,1}.
$$
$\phi$ extends to a homomorphism of stable groups from
$\beta_\infty := \lim_{g \to \infty} \beta _{2g}$ to
$\Gamma_\infty := \lim_{g \to \infty} \Gamma _{g,1}$. We study the
effect of this map in group (co)homology. \vskip .1in By work of
Fuks \cite{F1} and F. Cohen \cite{CLM} the (co)homology of the
braid groups has been well known. Though the (co)homology of the
finite genus mapping class groups remains difficult, through
recent work on a homotopy theoretic refinement of the Mumford
conjecture, the (co)homology of the stable mapping class group is
well understood rationally and integrally, cf. \cite{MT},
\cite{MW}, \cite{G}. In particular, both $\beta_\infty$ and
$\Gamma _\infty$ have ample integral (co)homology. \vskip .2in

\begin{thm}
The image of the map $\phi_*: H_* (\beta_\infty; \mathbb Z)
\rightarrow H_* (\Gamma_\infty; \mathbb Z)$ is zero for $*>0$.
\end{thm}

\vskip .1in This is our main result.  We believe  though that the
method of proof via categorical delooping is of interest in its
own right. In particular, we introduce and study the homotopy type
of a new monoidal category (operad) $\mathcal{ T}$, built out of
generalised Artin groups on trees and forests.
%Indeed, $\phi$ extends to these generalised Artin groups
%and the analogue of Theorem 1.1 holds.
\vskip .1in Most of the difficulty in proving Theorem 1.1 stems
from the fact that $\phi$ is defined  algebraically. In an
Appendix we prove an analogue of Theorem 1.1 for two different
geometrically defined maps. \vskip .3in

\noindent {\bf 1.2 Background: } In the early 1980s, in
conversation with E. Miller and F. Cohen, John Harer made the
following conjecture based on many explicit calculations he had
done using the Thurston-Hatcher method of giving a presentation of
the mapping class group. \vskip .3in

\noindent {\bf Conjecture (J. Harer)} {\it The map $\phi_*: H_*
(\beta_\infty ; \mathbb Z/2 \mathbb Z) \to H_*(\Gamma _\infty;
\mathbb Z /2 \mathbb Z)$ is trivial.} \vskip .2in The approach
suggested in \cite{C}, as ours in this paper, is based on product
structures. Juxtaposition of strings defines a monoidal product on
the disjoint union of classifying spaces of the braid groups,
$\coprod _{k\geq 0} B \beta_{k}$. Its group completion is the
space
$$
\mathbb Z \times B\beta _\infty ^+ \simeq \Omega ^2 S^2;
$$
here \lq\lq + " denotes Quillen's plus construction. The loop
space structure on $\Omega ^2 S^2$ is induced by  this monoidal
product. The additional, two-fold loop space structure
corresponds to a wreath product map of braid groups,
$$
\omega_\beta: \beta _q \wr \beta _k \longrightarrow \beta _{qk}
$$
which maps the  element $(\sigma; \mu_1, \dots , \mu_q)$ to the
braid on $qk$ strings where the braids $\mu_1, \dots , \mu _q$ are
weaved together (each considered a single strand) as prescribed by
$\sigma$. Similarly, a monoidal product can be defined on
$\coprod_{g\geq 0} B\Gamma _{g,1}$ by joining two surfaces by a
pair of pants surfaces. This induces a loop space structure on the
group completion, $\mathbb Z \times B \Gamma ^+_\infty$.  Miller
\cite{Mi} (see also \cite{Bo} and \cite{So})
noted that this extends to a double loop space structure
which corresponds to the wreath product
$$
\omega_\Gamma: \beta _q \wr \Gamma _{g,1} \longrightarrow \Gamma
_{qg, 1}.
$$
To define $\omega _\Gamma$ one identifies $\beta _q$ as a subgroup
of the mapping class group of a $q$-legged pants surface (i.e. a
sphere with $q+1$ disks removed) in which the $q$ legs can be
permuted, cf. \cite{CT}. \vskip .1in In homology the wreath
products give rise to the first Araki-Kudo-Dyer-Lashof operation
$Q_1: H_* \to H_{2*+1}$, and F. Cohen in \cite{CLM} has computed
$Q_1$ in the homology of the braid groups with $\mathbb Z
/q\mathbb Z$ -coefficients. In particular, he proves
$$
H_* ( \beta _\infty ; \mathbb Z /2\mathbb Z) \simeq \mathbb Z /
2\mathbb Z [ a_1, a_2, \dots]
$$
as polynomial algebras, where $| a_i| = 2^i -1$ and $Q_1 (a_i) =
a_{i+1}$. As by \cite{P} the first homology group of $\Gamma
_{g,1}$ is trivial for $g\geq 3$, Harer's conjecture would follow
if $\phi$ commutes with the wreath products (for $q=2$). However,
this is not the case. Moreover, Maginnis \cite{Mag} shows that
this is also not the case in homology: the difference between
$\phi_* \circ Q_1$ and $ Q_1 \circ \phi_*$, when applied to the
non-trivial 1-dimensional class, is a non-zero element in $H_3
(\Gamma _{2,1}; \mathbb Z / 2 \mathbb Z)$. \vskip .2in

\noindent {\bf 1.3 Strategy:} The fact that $\phi$ does not
commute with the wreath products means that $\phi : \coprod
_{g\geq 0} B\beta_{2g} \to \coprod _{g\geq 0} B\Gamma _{g,1}$ is
not a map of $E_2$-spaces in the sense of May  \cite{CLM}.
Nevertheless, we will show that the induced map on group
completion is. \vskip .2in \noindent {\bf Theorem 5.2} \/ \/ {\it
$\phi : B \beta_\infty ^+
 \longrightarrow B\Gamma ^+_\infty$ is a map of double loop
 spaces.}
\vskip .1in Theorem 1.1 and Harer's conjecture follow immediately;
see Lemma 5.3. We note here that Maginnis' obstruction lies in the
unstable range of the homology of the mapping class group, cf.
\cite{H}, \cite{I}, and therefore there is no contradiction.
\vskip .1in The idea of the proof is to use categorical delooping.
Recall, when a connected category $\mathcal {C}$ is monoidal
(symmetric) then its classifying space $B \mathcal {C}$ is a loop
space (infinite loop space, cf. \cite{Ma}, \cite{S2}).
Furthermore, monoidal (symmetric) functors induce maps  of loop
spaces (infinite loop spaces). By \cite{T1}, there is a symmetric
surface category $\mathcal {S}$ with $\Omega B\mathcal {S} \simeq
\mathbb Z\times B\Gamma ^+_\infty$. Our aim is therefore to
construct a monoidal category $\mathcal {T}$ out of braid groups
and extend $\phi$ to a monoidal functor $\Phi : \mathcal {T} \to
\mathcal {S}$; the induced map $\Omega B\Phi : \Omega B\mathcal
{T} \to \Omega B \mathcal {S}$ will automatically be a map of
double loop spaces. \vskip .2in

\noindent {\bf 1.4 The category $\mathcal {T}$:} The construction
of the category $\mathcal {T}$ imitates the construction of
$\mathcal {S}$. The place of surfaces in $\mathcal {S}$ is taken
by trees in $\mathcal {T}$. We are lead to consider generalised
Artin groups associated with any tree and forest. Recall, given a
graph $\Sigma $, its Artin group $ \beta (\Sigma )$ is the group
generated by the edges in $\Sigma$. Two generators $e,f$ commute
if they are disjoint, and satisfy the braid relation $efe=fef$ if
they have a vertex in common. The composition (operad structure)
in $\mathcal {T}$ is given by grafting of trees. The homomorphism
$\phi$ can now be extended to a functor from $\mathcal {T}$ to
$\mathcal S$ through a careful choice of images for the generators
at branch points; see Figure 7. Some complexity in the proof
arises from the fact that $\Omega B \mathcal {T}$ is not  homotopy
equivalent to $\mathbb Z\times B\beta^+_\infty$. For comparision,
we recall that Fuks \cite{F2} proves a homotopy equivalence
$$
B \beta (D_\infty) ^+ \simeq \Omega ^2 S^3 \times S^3 \{ 2\}
$$
where $S^3 \{ 2 \}$ denotes the homotopy fiber of the degree 2 map
from $S^3 $ to $S^3$ and $D_\infty$  the infinite Coxeter graph of
type $D$. To prove Theorem 5.3 it is however enough to show that
$\phi$ factors as a map of  double loop spaces  through the
identity component $\Omega _0 B\mathcal {T}$.  Indeed, in Section
4,  we prove \vskip .2in \noindent {\bf Theorem 4.1}\/ \/ {\it
$\Omega_0 B \mathcal {T} \simeq \Omega ^2 S^3 \times  \Omega ^2 W$
as double loop spaces for some $W$.} \vskip .2in The homotopy type
of $W$ remains undetermined.   \vskip .3in

\section { Definition of tile categories}
\vskip .2in The categories we will consider in this paper are most
naturally seen as  strict 2-categories. Recall, a {\it 2-category}
is a category $\mathcal {C}$ enriched over $CAT$, the category of
small categories. This means that instead of morphism sets,
$\mathcal {C}$ has morphism categories $\mathcal {C}(a,b)$ and
composition $\mathcal {C}(a,b) \times \mathcal {C} (b, c) \to
\mathcal {C}(a,c)$ is a functor for all objects $a,b,c, \in
\mathcal {C}$. A 2-category $\mathcal {C}$ is {\it strict} if
composition is strictly associative. We may apply the nerve
construction to the morphism categories. Since the nerve
construction commutes with products, this gives a category
$\mathcal {B}\mathcal {C}$ enriched over $\Delta$-sets, the
category of simplicial sets. In other words, $\mathcal {B}\mathcal
{C}$ is a simplicial category with constant objects. Applying the
nerve construction to $\mathcal {B}\mathcal {C}$ yields a
bisimplicial set. The realization of this bisimplicial set is the
{\it classifying space} $B\mathcal {C}$ of $\mathcal {C}$. \vskip
.1in The 2-category $\mathcal {T}$ may be thought of as a special
type of cobordism category. Its objects are the natural numbers
$\Bbb N$ where $n$ represents $n$ ordered copies of the unit
interval. Its 1-morphisms $\mathcal {T}(n,m)$ are generated by
three atomic tiles $D: 0\to 1, P: 2 \to 1$ and $F: 1 \to 1$, as
illustrated in Figure 2.
\begin{center}
\epsfig{figure=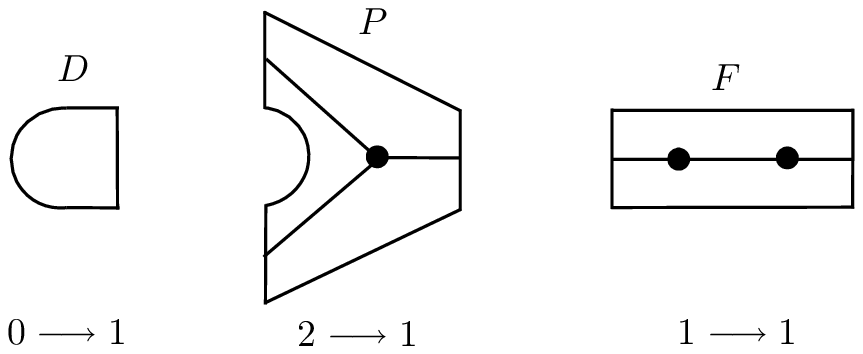, width=0.7\linewidth}

Figure 2.
\end{center}
\vskip .2in These are three discs thought of as cobordisms between
intervals. In addition, the tiles $P$ and $F$ have one
respectively two marked points. The 1-morphisms are generated by
these three atomic tiles by gluing along incoming and outgoing
intervals and disjoint union. Note that all incoming and outgoing
intervals are  ordered. A typical 1-morphism is illustrated in
Figure 3. Identity morphisms are added and may be thought of as
the objects themselves, i.e. zero length cobordisms. Composition
of 1-morphisms is defined by gluing. {\it Note:} The 1-morphisms
are  freely generated by the atomic tiles as morphisms of a
monoidal category. Given any two tiles $T_1, T_2$, we identify
$(T_1 \sqcup 1) \circ (1 \sqcup T_2) = (T_1 \sqcup T_2) = (1
\sqcup T_2) \circ (T_1 \sqcup 1)$. However, homeomorphic tiles are
not identified in general and disjoint union is {\it not}
symmetric. Indeed, $T_1 \sqcup T_2 = T_2 \sqcup T_1$ if and only
if $T_1 =T_2$. \vskip .1in The morphism categories $\mathcal {T}
(n,m)$ are disjoint unions of  groups. Thus there are no morphisms
between two different tiles, i.e. $\mathcal {T} (T_1, T_2) =
\emptyset$ if $T_1 \neq T_2$. We will now describe a functorial
choice of generators for the 2-endomorphism groups. Given a tile
$T \in \mathcal {T} (n,m)$, join the marked points to form a graph
$\Sigma_T$, as indicated in Figure 3.
\begin{center}
\epsfig{figure=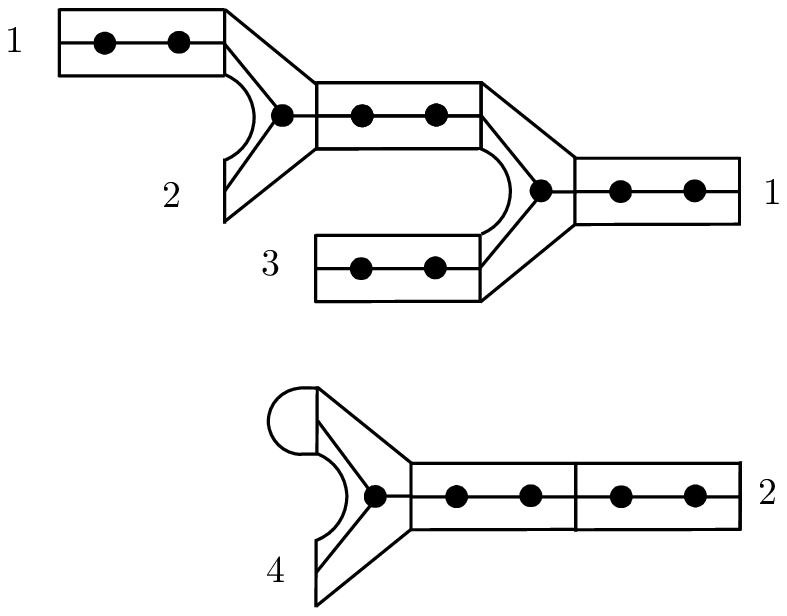, width=0.7\linewidth}

Figure 3.
\end{center}
\vskip .1in
 $\Sigma _T$ is  a union of $m$ disjoint trees and can be built
functorially  out of trees on atomic tiles, see Figure 2. Any
half-edges in $T$ that are not completed through gluing are not
considered to be part of $\Sigma _T$. Then by definition, the
construction of the graphs is functorial with respect to gluing
and disjoint union, i.e. for two tiles $T_1$ and $T_2$, $ \Sigma
_{T_1} $ and $\Sigma _{T_2} $ can canonically be identified as
subgraphs of $\Sigma _{T_1 \circ T_2}$ and $\Sigma _{T_1 \sqcup
T_2}$. Define $\mathcal {T} (T,T)$ to be the group generated by
the  edges $e, f \in \Sigma _T$ satisfying the following
relations.
$$
\aligned e f
= fe\quad &\text { if } \quad e\cap f = \emptyset; \\
efe=fef \quad &\text { if } \quad  e\cap f \neq \emptyset.
\endaligned
\leqno (2.1)
$$
Thus  $\mathcal {T}(T, T) = \beta (\Sigma_T)$ is the Artin group
associated with the graph $\Sigma _T$; for further discussion see
for example \cite{V}. In particular, the group associated with the
tile $F^k = F\circ \dots \circ F $  ($k$ times) is the usual braid
group on $2k$ strings,
$$
\mathcal {T} (F^k, F^k) = \beta (\Sigma_{F^k}) = \beta _{2k}.
$$
This completes the definition of $\mathcal {T}$. \vskip .2in We
will need some information about the homotopy type of $\mathcal
{T}$. For this purpose, we introduce an auxiliary tile category
$\hat {\mathcal {T}}$ whose homotopy type we are able to
determine, see Proposition 3.2 below. The objects  and 1-morphisms
of $\hat {\mathcal {T}}$ are the same as in $\mathcal {T}$. The
morphism categories are groupoids: For tiles $T_1$ and $  T_2$ in
$\mathcal {T}(n,m)= \hat {\mathcal {T}}(n,m)$,  the set of
2-morphisms $\hat{\mathcal {T}}( T_1, T_2)$ is the set of
connected components of the space of homeomorphisms from $T_1$ to
$T_2$ that identify the ordered incoming and outgoing boundary
intervals, and map marked points bijectively to marked points.
(The edges of the associated trees may not be  preserved.) Hence, if
$T_1$ and $T_2 $ are not homeomorphic or don't have the same
number of marked points, $\hat {\mathcal {T}} (T_1, T_2) =
\emptyset$. On the other hand, the endomorphisms of a tile $T$
with $k$ marked points is the mapping class group of $T$:
$$
\hat{\mathcal {T}}(T, T) = \Gamma (T) = \beta _k.
$$
\vskip .1in In order to be able to compare the homotopy types of
the two categories we define a functor
$$
\Theta : \mathcal {T} \longrightarrow \hat{\mathcal {T}}.
$$
On objects and 1-morphisms $\Theta$ is the identity, and we are
left to define a homomorphism of 2-morphism groups associated to a
tile $T$. For each generator  $e \in \mathcal {T}(T,T)= \beta
(\Sigma _T)$, let  $\sigma _e$ in $\hat {\mathcal {T}} (T,T) =
\Gamma (T)$ be the isotopy class of a half Dehn twist, the
homeomorphism of $T$ that turns the edge $e$ by a half rotation in
a clockwise direction and is the identity outside a contractible
neighborhood of $e \subset T$. By definition this assignment is
functorial under gluing. It is also not difficult to see that the
assignment $e \mapsto \sigma_e$ induces a well-defined surjective
group homomorphism
$$
\Theta: \beta (\Sigma_T) \lra \Gamma (T).
$$
We expect this to be well-known and only sketch an argument. All
trees that occur have vertices that are at most trivalent. By an
inductive argument, it is enough to consider the tree $T$ in
Figure 4 with five vertices, $ 1,2,3,4,5 $, and edges $(1,2),
(2,4) (3,4), (4, 5)$. The corresponding generators are mapped  to
four elements $\sigma _1, \tilde \sigma, \sigma_3, \sigma_4$ of
$\Gamma (T) \simeq \beta _5$.
\begin{center}
\epsfig{figure=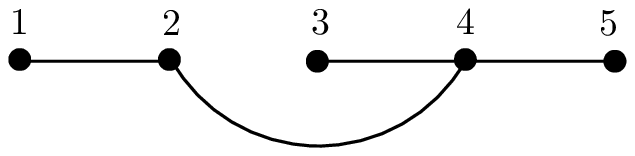, width=0.4\linewidth}

Figure 4.
\end{center}
\vskip .1in \noindent $\tilde \sigma $ corresponds in the standard
representation of $\beta _5$ to $\sigma _3 \sigma _2 \sigma
_3^{-1}$. It remains to check that the three relations in $\beta
(T)$ are satisfied in $\beta _5$. For example:
$$
\aligned \tilde \sigma \sigma _3 \tilde \sigma &= (\sigma_3 \sigma
_2 \sigma _3 ^{-1}) \sigma _3 (\sigma _3 \sigma _2
\sigma _3^{-1}) \\
&= \sigma _3 (\sigma _2 \sigma _3 \sigma _2) \sigma _3 ^{-1}    \\
&= \sigma _3 ^2 \sigma _2   \\
&= \sigma _3 (\sigma _3 \sigma _2 \sigma _3 ^{-1}) \sigma _3    \\
&= \sigma _3 \tilde \sigma \sigma _3.
\endaligned
$$
\vskip .2in \noindent {\bf Remark 2.2} The homomorphism $\Theta$
is in general  {\it not} injective. To invert the maps for the
tree in Figure 4, we would need that the element corresponding to
$\sigma _2= \sigma _3 ^{-1} \tilde \sigma \sigma _3$  commutes
with the one corresponding to $\sigma_4$. This, however, cannot be
deduced from the relations (2.1). Note that this is the only
additional relation satisfied by these generators in $\beta_5$.
\vskip .2in Both categories $\CalT$ and $\hat \CalT$ contain the
subcategory $\CalM$ with one object, $1$,  and 1-morphisms
generated by the atomic tile $F$:
$$
\CalM \subset \CalT \overset \Theta  \lra \hat \CalT.
$$
\vskip .1in Disjoint union of cobordisms makes $\CalT$ and  $\hat
\CalT$ into strict monoidal 2-categories. Furthermore, the functor
$\Theta: \CalT \to  \hat \CalT$ is a  map of strict monoidal
2-categories. Note that these monoidal structures are {\it not}
symmetric as there are no morphisms that permute the ordered
intervals and therefore no symmetries exist. Also, the subcategory
$\CalM$ is {\it not} monoidal. \vskip .4in

\section{ Classifying spaces of tile categories}
\vskip .1in We will now study the homotopy type of these
categories. The realization of $\CalB \CalM$ is the monoid
$\coprod _{k\geq 0} B \beta_{2k}$. Hence, as is well known, by the
group completion theorem for monoids, we have \vskip .2in
\begin{prop} $
\Omega B\CalM \simeq 2\mathbb Z \times B\beta^+_\infty. $
\end{prop}
\vskip .1in In order to calculate the homotopy type of $\hat
\CalT$ a generalisation of the group completion theorem for
categories will be used, cf. \cite{T1}, \cite{Mo}. The main idea
goes back to McDuff and Segal's approach to the group completion
theorem via homology fibrations \cite{MS}. \vskip .2in
\begin{prop} $\Omega B \hat \CalT \simeq \mathbb Z \times B
\beta_\infty ^+$.
\end{prop}
\begin{proof} $\CalB \hat \CalT (n,1)$ is the nerve of a groupoid
in which each connected component is determined by the number of
marked points of the underlying tiles; note, all such tiles are
connected. Hence, its homotopy type is given by $\coprod _{k\geq
n-1} B \beta_{k}$. Gluing the atomic tile $F \in \hat \CalT (1,1)$
on the right defines a self map  of $\CalB \hat \CalT (n,1)$. Let
$$
\CalB \hat \CalT_\infty (n) = \underset {F}  {\rm hocolim}\,\,
\CalB \hat \CalT(n,1) \simeq \mathbb Z \times B\beta _\infty
$$
be the corresponding homotopy colimit. Both $\CalB \hat \CalT (\_
, 1)$ and $\CalB \hat \CalT_\infty (\_)$ define contravariant
functors from $\CalB \hat \CalT$ to $\Delta$-sets. The homotopy
colimit of the first is the nerve of Quillen's  category of
objects over 1 in $\CalB \hat \CalT$,  and is therefore
contractible as this category has a final element. Hence,
$$
\underset {\CalB \hat \CalT}  {\rm hocolim}\,\, \CalB \hat \CalT
(\_ , 1) \simeq * .
$$
But then also the homotopy colimit of the second functor is
contractible as the two homotopy colimits commute and the nerve of
the translation category is contractible:
$$
\underset {\CalB \hat \CalT}  {\rm hocolim}\,\,  \CalB \hat
\CalT_\infty (\_ ) = \underset {F}  {\rm hocolim}\,\, (\underset
{\CalB \hat \CalT}  {\rm hocolim}\,\, \CalB \hat \CalT (\_ , 1) )
\simeq *.
$$
Next consider the canonical forgetful map
$$
\underset {\CalB \hat \CalT}  {\rm hocolim}\,\,  \CalB \hat
\CalT_\infty (\_ ) \overset {\pi} { \lra} B \hat \CalT
$$
to the nerve of $\CalB \hat \CalT$. Each 1-morphism in $\hat \CalT
(n,m)$ acts as a homology equivalence on the fiber
$$
\CalB \hat \CalT _\infty (m) \simeq \mathbb Z \times
B\beta_\infty.
$$
This is immediate as the homotopy type of the fiber does not
depend on the object $m$ and left translation is conjugate to
right translation in the braid group. Thus, $\pi$ is a homology
fibration, i.e. the canonical map of the fiber $\CalB \hat
\CalT_\infty (m)$ into the homotopy fiber ${\rm hofib} (\pi)
\simeq \Omega B \hat \CalT$ is a homology isomorphism. By the
Whitehead theorem for simple spaces this  gives the desired
homotopy equivalence after plus construction.  \qed
\end{proof}

\vskip .2in \begin {cor} The inclusion $\CalM \subset \hat \CalT$
induces the natural inclusion
$$
\Omega B \CalM \simeq 2\mathbb Z \times B\beta^+_\infty
\hookrightarrow \Omega B \hat \CalT \simeq \mathbb Z \times
B\beta^+_\infty.
$$
\end{cor}
\vskip .2in
\begin{proof} Define $\CalB \CalM_\infty (\_)$ in
analogue to $\CalB\hat \CalT _\infty (\_)$ in the proof above. The
inclusion $\CalM \subset \CalT$ then induces a map of homology
fibrations
$$
\underset {\CalB\CalM}  {\rm hocolim}\,\,  \CalB \hat \CalM_\infty
(\_ ) \lra \underset {\CalB \hat \CalT}  {\rm hocolim}\,\,  \CalB
\hat \CalT_\infty (\_ )
$$
which is a homotopy equivalence of total spaces  (as both are
contractible) and is homotopic to the inclusion $2\mathbb Z \times
B\beta _\infty \to \mathbb Z \times B\beta _\infty$ on the fiber
over the object 1. The corollary follows as before by an
application of Whitehead's theorem for simple spaces. \qed
\end{proof}
\vskip.3in

\section{ Double loop space structures} \vskip .1in

The monoidal structures on $\CalT$ and $\hat \CalT$ induce double
loop space structures on $\Omega B \CalT$ and $\Omega B\hat
\CalT$. A priori this double loop space structure may however not
be compatible with  the standard double loop space structure on
$$
B\beta^+_\infty \simeq \Omega_0^2 S^2 \simeq \Omega ^2 S^3;
$$
here the subscript $0$ indicates that only the $0$-component is
considered. Below $\Omega ^2 S^3$ will be considered with its
standard double loop space structure. \vskip .2in
\begin {thm} There is a splitting of double loop
spaces
$$
\Omega _0 B \CalT \simeq \Omega ^2 S^3 \times \Omega ^2 W
$$
for some 2-connected space $W$. Furthermore, the inclusion $\CalM
\subset \CalT$ induces a map of double loop spaces
$$
\Omega ^2 S^3 \simeq \Omega _0 B\CalM \lra \Omega B \CalT
$$
which is homotopic to the inclusion of the first factor in the
above splitting.
\end{thm}
\vskip .2in \begin{proof} We first recall that
$$
\Omega _0 B \CalM \simeq \Omega ^2 S^3
$$
is a homotopy equivalence of loop spaces, cf. \cite{S1},
\cite{CLM}. Consider now the maps of $0$-components
$$
\Omega _0 B\CalM \lra \Omega_0 B \CalT \overset {\Omega B \Theta}
 {\lra} \Omega_0 B \hat \CalT.
$$
The first is a map of loop spaces by definition, while the second
is a map of double loop spaces as $\Theta $ is monoidal. Thus, by
Corollary 3.3, the composition is a homotopy equivalence of loop
spaces. By Lemma 4.3 below, the loop space structure on $\Omega
S^3$ is unique. Hence,
$$
\Omega _0 B \hat \CalT \simeq \Omega ^2 S^3
$$
is a homotopy equivalence of double loop spaces. \vskip .1in As
$\Omega B \Theta$ is a map of double loop spaces, its homotopy
fiber is a connected double loop space $\Omega ^2 W$ for some
2-connected space $W$. Because  $\Omega_0 B \CalT$ is an H-space,
the retraction above gives rise to a splitting of spaces,
$$
\Omega _0 B \CalT \simeq \Omega ^2 S^3  \times \Omega ^2 W; \leqno
(4.2)
$$ here the projection onto the first factor and the inclusion of
the second factor are maps of double loop spaces. \vskip .1in
$S^1$ embeds canonically in $\Omega ^2 S^3$. Let $f: S^1 \to
\Omega _0 B \CalT$ be the restriction of the inclusion of the
first factor in (4.2) to $S^1$ with respect to the splitting
(4.2), $f$ may be written as $(f_1, f_2)$ where  $f_1$ is the
inclusion and $f_2$ is the trivial map. $\Omega ^2 S^3 = \Omega ^2
\Sigma ^2 (S^1)$ is the free double loop space on $S^1$. By the
universal property for free double loop spaces, $f$ extends (up to
homotopy) uniquely to a map of double loop spaces  $\bar f :
\Omega ^2 S^3 \to \Omega _0 B \CalT$. Let $\bar f = (\bar f_1 ,
\bar f_2)$ with respect to the splitting (4.2). $\Omega ^2 W$ is a
double loop subspace of $\Omega_0 B \CalT$, and hence $\bar f_2$
is a map of double loop spaces extending $f_2 \simeq *$. Thus by
the uniqueness property, $\bar f_2$ is therefore trivial as well.
Finally, $\bar f_1 = \Omega B\Theta \circ \bar f$ is a map of
double loop spaces of  $\Omega ^2 S^3$ to itself. As it extends
$f_1$, by the uniqueness property, $\bar f_1$ is homotopic to the
identity. Hence, the splitting (4.2) is a splitting of double loop
spaces. In particular the inclusion of the first factor, which was
induced by the inclusion $\CalM \subset \CalT$, is  a map of
double loop spaces. \qed
\end{proof}
\vskip .2in \noindent {\bf Lemma 4.3}\,\,\, {\it For $n$ odd, the
loop space structure on $\Omega S^n$ is unique up to homotopy.}
\vskip .1in \begin{proof} For $n$ odd, the cohomology of $\Omega
S^n$ is a divided power algebra, i.e. there are generators $z_k\in
H^{k(n-1)} (\Omega S^n)= \mathbb Z$ such that
$$
k! \, z_k = z_1 ^k.
$$
Let $Z = \Omega S^n$ with some (potentially exotic) loop space
structure, and let $f: S^{n-1} \to Z$ be a generator of $\pi
_{n-1} Z = \mathbb Z$. By the universal property of free loop
spaces, $f$ extends to a map $\bar f : \Omega \Sigma (S^{n-1})
=\Omega S^n \to Z$ of loop spaces which is unique up to homotopy.
Clearly,
$$
\bar f^* (z_1) = f^* (z_1) = z_1,
$$
and,  as $\bar f^*$ is a map of cohomology rings,
$$
k! \,  \bar f^* (z_k) = \bar f^* (z_1 ^k) = z_1^k = k!\, z_k.
$$
Hence $\bar f^*(z_k) = z_k$ as $H^* (\Omega S^n)$ is torsion free.
Thus $\bar f$ is a homology equivalence and by Whitehead's theorem
for simple spaces a homotopy equivalence. It follows that $Z$ and
$\Omega S^n$ are homotopy equivalent as loop spaces.  \qed
\end{proof}
\vskip .2in \noindent {\bf Remark 4.4} The loop space structure on
$\Omega S^2$, and hence  the double loop space structure on
$\Omega ^2 S^2$, are not unique which can be seen as follows. The
generator of the first cohomology group defines a map $g: \Omega
S^2 \to S^1$ which is a splitting of the embedding $f: S^1 \to
\Omega S^2$. This induces a splitting of spaces
$$
\Omega S^2 \simeq S^1 \times F
$$
where $F$ is the homotopy fiber of $g$. $F$ has a natural loop
space structure as $g$ is a map of loop spaces. Indeed, $F\simeq
\Omega S^3$ as $g$ is induced by the Hopf map
$$
\Omega S^3 \lra \Omega S^2 \overset g  \lra S^1 \lra S^3 \lra S^2.
$$
However, the loop space structure of $\Omega S^2$ is not the same
as the product loop space structure of $S^1 \times \Omega S^3$,
for $S^2$ is not homotopy equivalent to $\Bbb CP^\infty \times
S^3$. 
\vskip .1in 
The same argument shows that more generally the $\Omega ^{n-1}$-structure on 
$\Omega ^{n-1} S^n$ and the $\Omega ^n$-structure of $\Omega ^n S^n$
are not unique as $S^n$ does not have  $K(\Bbb Z, n )$ as a retract.

\vskip .2in\noindent{\bf Remark 4.5} Let $\CalT _0$ denote the
subcategory of $\CalT$ with the same objects and 1-morphisms but
only identity 2-morphisms. As there are only endomorphisms in
$\CalT$, there is a natural retraction $\pi_0: \CalT \to \CalT_0$
mapping each 2-morphism group to its identity. $\CalT_0$ inherits
a monoidal structure from $\CalT$ and both the inclusion and the
projection are monoidal. Thus, for some space $Z$, there is a
splitting of double loop spaces
$$
\Omega B \CalT \simeq \Omega B \CalT_0 \times \Omega ^2 Z.
$$
Furthermore, from the proof of Theorem 4.1, it follows that for
some two connected space $\tilde W$,
$$
\Omega ^2 Z \simeq \Omega ^2 S^3 \times \Omega ^2 \tilde W.
$$
We expect that $\Omega B \CalT_0$ is homotopic to a discrete,
free abelian group on
infinitely many generators. This would, in particular, imply that
$W =\tilde W$. \vskip .4in
\section{ Surface category and the functor $\Phi$}
\vskip .1in We define a surface 2-category $\CalS$ following
\cite{T1} and \cite{T2} which gives rise to a convenient
categorical delooping of $B\Gamma ^+_\infty$. Its objects are the
natural numbers $\Bbb N$ with each $n$ representing $n$ disjoint,
ordered circles. The 1-morphisms are generated by three atomic
surfaces, a disk $D': 0\to 1$, a torus with two incoming and one
outgoing boundary component $P': 2\to 1$, and a torus with one
incoming and outgoing boundary component $T': 1 \to 1$.
\begin{center}
\epsfig{figure=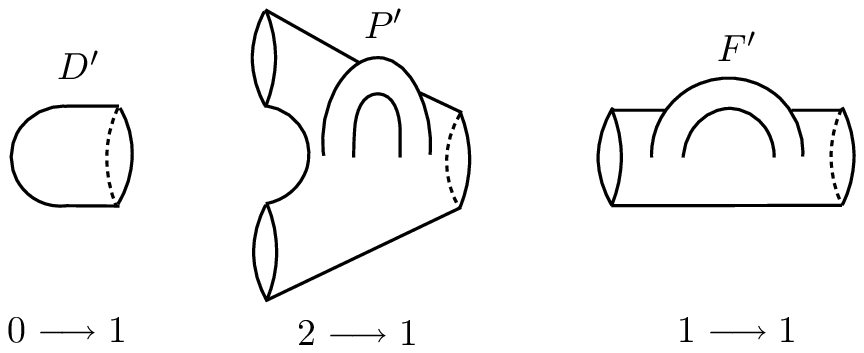, width=0.7\linewidth}

Figure 5.
\end{center}
\vskip .1in In analogy with the tile category $\CalT$, any
1-morphism is obtained from these atomic surfaces by gluing
incoming boundary circles of one surface to outgoing boundary
circles of another surface and by disjoint union. In addition
reordering of the incoming or outgoing boundary circles gives rise
to new 1-morphisms. Identity 1-morphisms are adjoined  and may be
thought of as zero length cylinders, i.e. the objects themselves.
A typical 1-morphism is illustrated in Figure 6.
\begin{center}
\epsfig{figure=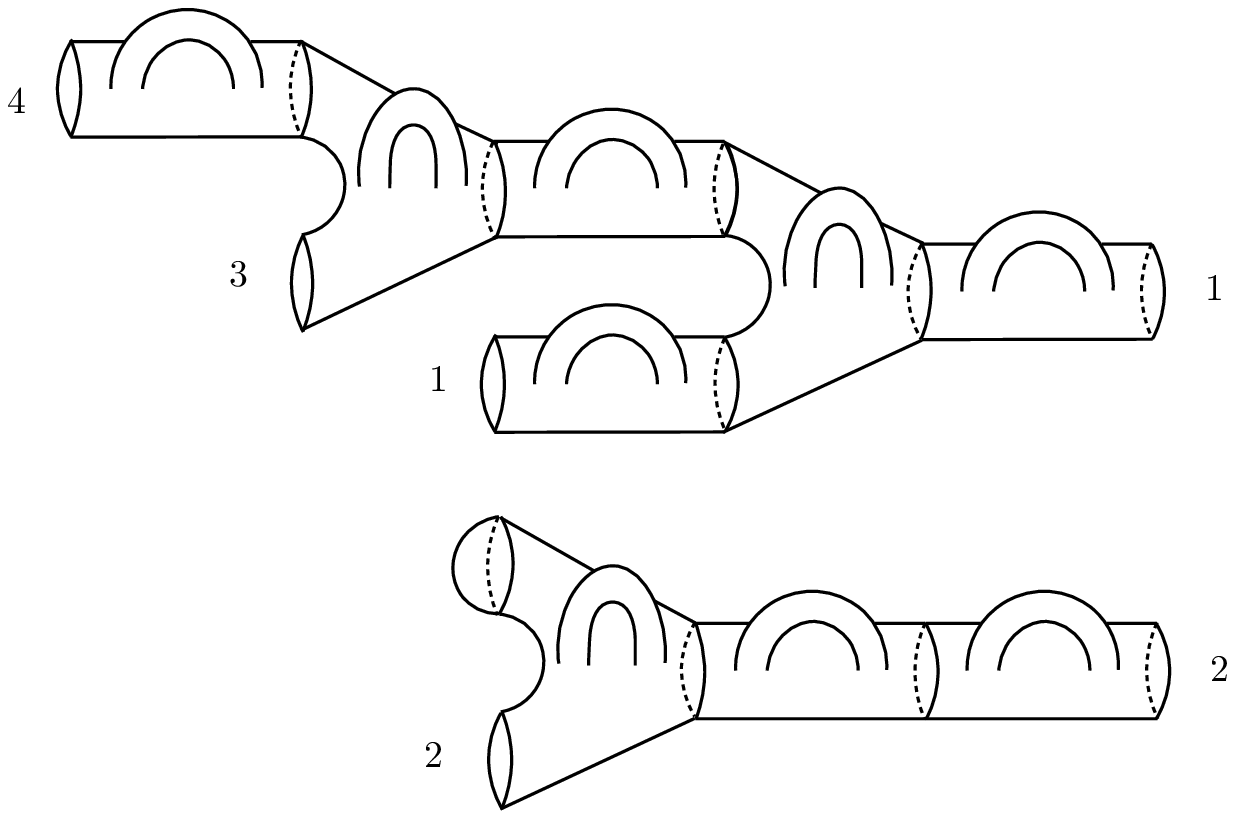, width=0.7\linewidth}

Figure 6.
\end{center}
\vskip .1in The set of 2-morphisms $\CalS( S_1, S_2)$ between two
1-morphisms $S_1, S_2 \in \CalS (n,m)$ is the set of connected
components of the space of homeomorphisms from $S_1$ to $S_2$
which identify the $n$ incoming and $m$ outgoing boundary circles.
When $S_1 =S_2$ this is the  mapping class group of the surface.
\vskip .1in Let $\CalM_\Gamma$ denote the subcategory of $\CalS$
with only on object, $1$, 1-morphisms generated by $ F'$ and a
full set of 2-morphisms. Recall from \cite{T1}, \cite{T2} that the
group completion argument in conjunction with Harer's homology
stability theorem \cite{H} for the mapping class groups proves
\vskip .2in
\begin{thm} $\Omega B\CalM_\Gamma \simeq \Omega B \CalS \simeq
\mathbb Z \times B\Gamma ^+_\infty.$
\end{thm}
\vskip .2in Disjoint union makes $\CalS$ into a strict symmetric
monoidal 2-category. As $\CalS$ is connected, its classifying
space $B\CalS$ and thus its loop space are therefore  infinite
loop spaces. Wahl \cite{Wa} proved that the induced infinite loop
space structure on $\mathbb Z\times B\Gamma ^+_\infty$ is
compatible with the double loop space structure defined  by the
pairs of pants product discussed  in Subsection 1.1 of  the
introduction. \vskip .2in We will now define the monoidal
2-functor $ \Phi: \CalT \to \CalS. $ On objects, $\Phi$ maps $n$
intervals to $n$ circles, i.e. $\Phi (n) = n$. For the atomic
tiles, define
$$
\Phi (D) = D', \quad \Phi (P) = P', \quad \Phi (F) = F'.
$$
$\Phi$ extends to a map of all 1-morphisms as by definition the
1-morphisms in $\CalT$ and $\CalS$ are built from the atomic
1-morphisms in the same fashion. Under this map Figure 6 is the
image of Figure 3 apart from a different labelling of the incoming
boundary components. \vskip .2in To define the functor $\Phi$ on
2-morphisms, recall from Section 2 that a tile $T$ has an
associated tree $\Sigma_T$. $\Phi$ will map this tree to a system
of curves on the surface $\Phi (T)$. For the atomic tiles $P$ and
$F$ this is indicated in Figure 7.
\begin{center}
\epsfig{figure=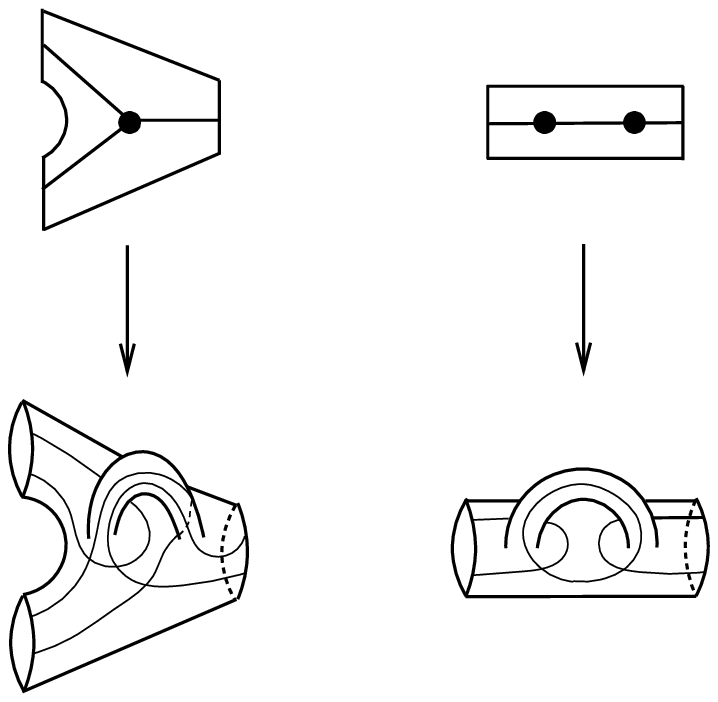, width=0.7\linewidth}

Figure 7.
\end{center}
\vskip .1in Through gluing this assignment extends to more general
tiles. Half edges and half circles gain importance when under
gluing they are completed and give rise to additional generators.
Otherwise, they should be ignored. Thus, each edge $e \in \Sigma
_T$ is  mapped to a simple closed curve $C_e$, and two such curves
$C_e$ and $C_f$ intersect (precisely once) if and only if the
edges $e$ and $f$ share a vertex. Dehn twists around curves that
don't intersect commute. It is also well-known that Dehn twists
around two simple closed curves which intersect  once satisfy the
second braid relation (2.1), cf. \cite{B}. Thus the map of
generators $e\mapsto C_e$ extends to a well-defined group
homomorphism
$$
\Phi : \CalT (T,T)= \beta (\Sigma _T)  \lra \CalS(\Phi (T), \Phi
(T)) = \Gamma (\Phi (T)).
$$
By definition the assignment of curve systems to trees is
functorial under gluing. Furthermore, $\Phi$ clearly commutes with
disjoint union, and hence defines a (strict) monoidal 2-functor.
\vskip .2in
%\noindent
%{\bf Remark 5.2.}
%We don't expect the Dehn twists to satisfy the
%additional relations satisfied by the
%edges of $T$ interpreted as elements
%in $\Gamma (T)$, cf. Remark 1.3.
%Furthermore,
%because of Maginnis' obstrcution [M], the functor
%$\Phi$ can {\it not} be extended to
%category $\hat \Cal T$.
%This promted the introduction of  the category $\Cal T$.
\begin{thm} The map $\phi : B\beta^+_\infty \to B\Gamma
^+_\infty$ is a map of double loop spaces.
\end{thm}
\vskip .2in \begin{proof} On the subcategory $\CalM$ of $\CalT$
the functor $\Phi$ is identical to the homomorphism $\phi:
\beta_{2g} \to \Gamma_{g, 1+1}$ of the introduction. Hence the
following diagram commutes up to homotopy.
$$
\CD 2\mathbb Z \times B\beta^+_\infty \simeq \Omega B\CalM @>>>
    \Omega B\CalT \\
@ V \phi VV @VV \Omega B \Phi V \\
\mathbb Z \times B \Gamma ^+_\infty \simeq \Omega B\CalM_\Gamma @>
\simeq >>
    \Omega B \CalS.
\endCD
$$
By Theorem 4.1, the restriction to the 0-component of the top
horizontal map is a map of double loop spaces. As $\Phi$ is
monoidal,  $\Omega B \Phi$ is a double loop space map, and hence
so is $\phi$ when restricted to the 0-component. \qed
\end{proof}
\vskip .2in \begin{lem} Any map $\varphi: B\beta ^+_\infty \to
B\Gamma ^+_\infty$ of double  loop spaces is homotopically
trivial.
\end{lem}
\vskip .2in \begin{proof} As $B\beta^+_\infty \simeq \Omega ^2
S^3$ is the free double loop space on the circle $S^1$, and as
$\varphi$ is a map of double loop spaces, up to homotopy it
is determined
completely by its restriction to $S^1$. However, by \cite{P},  the
mapping class group is perfect for $g\geq 3$, and $B\Gamma
^+_\infty$ is simply connected. Hence, the restriction of
$\varphi$ to $S^1$ is homotopically trivial, and so is $\varphi$.
\qed
\end{proof}
\vskip .1in As the plus construction does not change the homology,
this  also implies the triviality of the map $\varphi$ on group
homology. In particular, for $\varphi= \phi$ this proves Theorem
1.1 of the introduction and Harer's conjecture as a special case.
By \cite{H} and \cite{I}, we can state an   unstable  version of
this result. \vskip .2in \begin {cor} The image of $\phi:
H_*(\beta_{2g}; \mathbb Z) \to H_* (\Gamma _{g, 1};\mathbb Z)$ is
zero for $0< * < g/2$, and the image of any element of degree
$*\geq g/2$ is zero or unstable.
\end{cor}
%\vskip .1in
%We will now show that the map $\Omega B \Phi$ is null homotopic
%on all of $\Omega_0 B \CalT$. Let $\bar \CalT$ be the subcategory
%of $\hat \CalT$ with the same objects and 1-morphisms but where the
%2-morphisms are restricted to endomorphisms.
\vskip .4in

\section{ Appendix: Geometrically defined homomorphisms}
\vskip .1in In this appendix we discuss  maps from the braid group
to the mapping class group which are defined geometrically, i.e.
by identifying the braid group as (a subgroup of) the mapping
class group of a subsurface. There are many such maps but we will
discuss only two in detail. Others can be analyzed in  a similar
fashion. \vskip .1in The basic idea is to identify the braid group
as a subgroup of the mapping class group of a genus zero surface
with boundary components as follows. Let $S_{0, k+1}$ be a sphere
with $k+1$ disks removed and parametrised boundary circles
$\partial _0, , \partial _1, \dots, \partial _k$. Consider the
orientation preserving diffeomorphisms that fix the first boundary
component $\partial _0$ pointwise but may permute the other $k$
boundary components as long as they preserve the parametrisation
of each. The associated mapping class group, denoted by $\Gamma
_{0, (k), 1}$, is the ribbon braid group $R\beta _k$ on $k$
ribbons. $R\beta _k$ is the wreath product $\beta _k \wr \mathbb
Z$, and $\beta _k$ can naturally be identified as a subgroup.
$$
\beta_k \subset \beta_k \wr \mathbb Z = R\beta_k \simeq \Gamma
_{0, (k),1}.
$$
\vskip .1in There are  two ways that this identification leads to
homomorphisms of the braid group into the mapping class group that
we consider here. The first is obtained by gluing  $k$ copies of
the surface $S_{g,1}$ onto $S_{0, k+1}$ along $\partial_1, \dots,
\partial_k$, and then extending the diffeomorphism by the identity
and permutation the $k$ copies of $S_{g,1}$, cf. Figure 8. \vskip
.1in
\begin{center}
\epsfig{figure=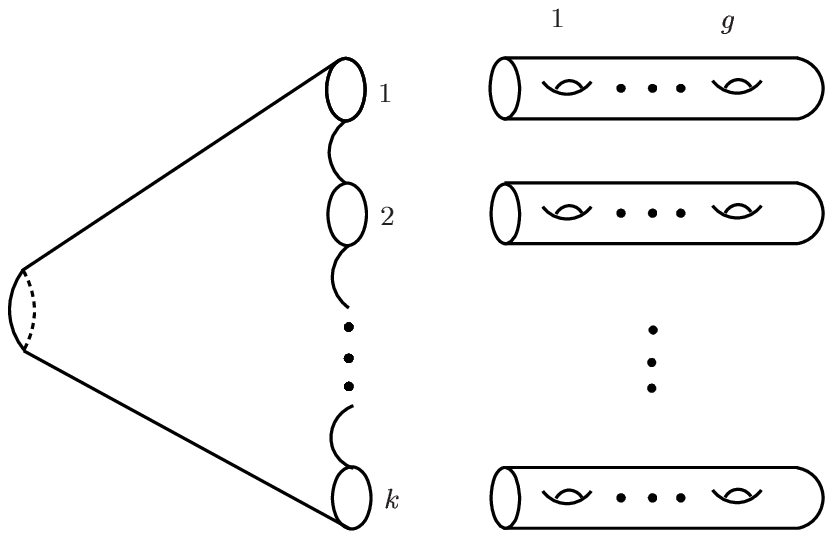, width=0.7\linewidth}

Figure 8.
\end{center}
\vskip .1in On mapping class groups this map factors through the
wreath product map $\omega_\Gamma$ described in the introduction
$$
\varphi_1: \beta_k \hookrightarrow \beta_k \wr \Gamma _{g,1}
\overset {\omega_{\Gamma}} \lra \Gamma _{kg, 1}.
$$
As these give rise to Miller's double loop space structure
\cite{Mi}, the induced map on group completions is a map of double
loop spaces:
$$
\varphi_1: \Omega ^2 S^2 \lra \Omega^\infty S^\infty \lra \mathbb
Z \times B\Gamma ^+_\infty.
$$
It follows from Lemma 5.3 that \vskip .2in
\begin{prop}
The image of $\varphi_1: H_*(\beta _k; \mathbb Z) \to H_* (\Gamma
_{kg, 1}; \mathbb Z)$ is zero for $0<* < kg/2$.
\end{prop}
\vskip .1in The second family of homomorphisms is constructed as
follows. Consider two copies of the surface $S_{0, k+1}$ glued
along their boundary components $\partial _1, \dots , \partial _k$
to form a surface $S_{k-1, 2}$ which in turn is a subsurface of a
larger surface $S_{g +k,n}$. We will only consider the case when
one of the two boundary components $\partial _0$ remains a
boundary component in $S_{g+k ,n}$ with $n= 2$, cf. Figure 9, and
leave the other cases as an exercise. \vskip .1in
\begin{center}
\epsfig{figure=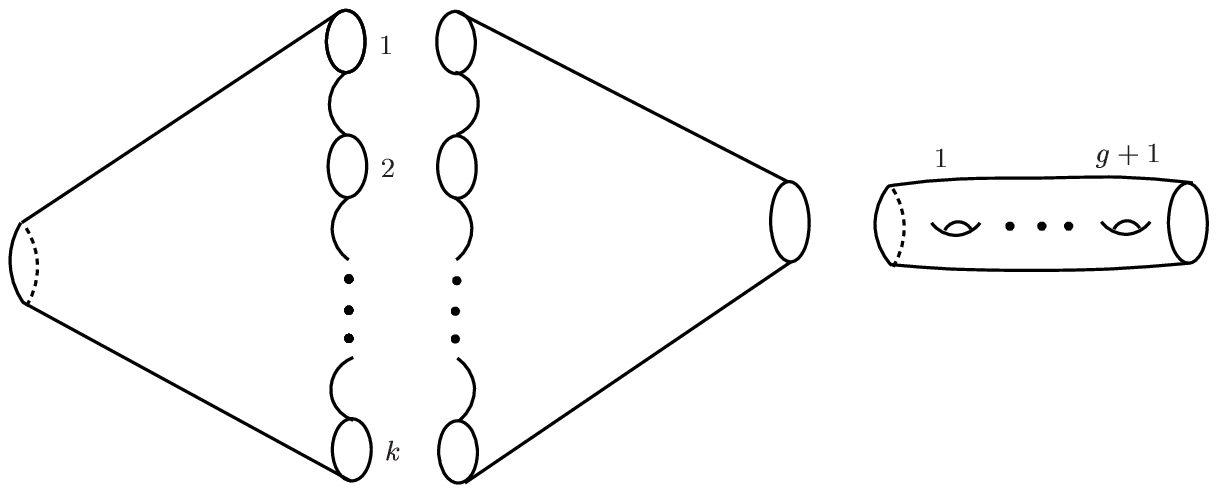, width=0.9\linewidth}

Figure 9.
\end{center}
\vskip .1in Any diffeomorphism of $S_{0,k}$ as described above can
be extended to $S_{k-1, 2}$ by \lq\lq mirroring" the action on the
second copy of $S_{0,k}$, and can then be extended to $S_{g +k,
n}$ by the identity diffeomorphism. This gives rise to a group
homomorphism
$$
\varphi_2: \beta_k \lra \beta _k \times_{\Sigma _k} \beta _k  \lra
\Gamma _{ g+k, 2},
$$
where the  group in the middle is defined as the pull-back in the
following diagram:
$$
\CD
\beta_k \times _{\Sigma_k} \beta_k  @>>>    \Sigma _k   \\
@VVV    @V\triangle VV  \\
\beta_k \times \beta_k      @>\pi \times \pi >>
        \Sigma _k \times \Sigma _k.
\endCD
$$
Here $\triangle$ denotes the diagonal map, and $\pi : \beta_k \to
\Sigma _k$ is the canonical surjection. \vskip .2in
\begin{prop}
The image of $\varphi_2: H_* (\beta _k; \mathbb Z) \to H_* (
\Gamma _{g+k,2}; \mathbb Z)$ contains at most 2-torsion for $0< *
< (g+k)/2$.
\end{prop}
\vskip .1in
\begin{proof}
Consider the following commutative diagram of group homomorphisms
$$
\CD
\beta _k    @>>> \beta_k \times_{\Sigma _k} \beta_k @>>> \Gamma _{g+k,2}    \\
@V s VV @V s\times s VV @V s VV \\
\Gamma _{g+k,  (k), 1} @>>>  \Gamma _{g+k, (k), 1} \times _{\Sigma
_k} \Gamma _{g+k,(k), 1}
    @>>> \Gamma _{2(g+k) +(k-1), 2}.
\endCD
$$
Here $\Gamma _{g+k, (k), 1}$ denotes the mapping class group of
the surface $S_{g+k, k+1}$ where   $k$ of the boundary components
may be permuted as long as the parametrisation of each component
is preserved while one of the boundary components is fixed
pointwise; the group in the middle on the bottom is defined as a
pull-back as above. The map $s$  is induced by gluing one of the
boundary components of $S_{g+k,2}$ along $\partial _0$ to  a copy
of $S_{0, k+1}$. The left  horizontal maps are defined by \lq \lq
mirroring", while the bottom right horizontal map is defined by
identifying the two copies of the boundary components $\partial
_1, \dots, \partial _k$. As a consequence of Harer-Ivanov
stability \cite{H},\cite{I}, it was proved in \cite{BT} that the
map $\beta _k \to \Gamma _{g+k, (k), 1}$ factors  in homology in
degrees $*< (g+k)/2$ through $\Sigma_k$. The right most vertical
map is a homology isomorphism in these degrees, and hence the
claim follows as the image of $\pi: \beta _k \to \Sigma _k$ in
homology contains only 2-torsion, cf. \cite{CLM}. \qed
\end{proof}
\vskip .2in We expect that the image is actually zero in the given
range though a proof of this seems at the moment out of reach.
\vskip .4in

\begin{acknowledgements} The second author
would like to thank Marcel B\"okstedt for his interest and
comments. \end{acknowledgements}
% ----------------------------------------------------------------

\end{document}